\documentclass[twoside,a4paper,12pt]{article}
\usepackage{amsmath,amsthm,amsfonts,amssymb}
\usepackage{mathrsfs}
\usepackage{color,graphicx}

\makeatletter
\def\@seccntformat#1{\csname the#1\endcsname.\quad}
\renewcommand\section{\@startsection {section}{1}{\z@}%
                                   {-3.5ex \@plus -1ex \@minus -.2ex}%
                                   {2.3ex \@plus.2ex}%
                                   {\normalfont\bf\center }}
\renewcommand\subsection{\@startsection {subsection}{1}{\z@}%
                                   {-3.5ex \@plus -1ex \@minus -.2ex}%
                                   {2.3ex \@plus.2ex}%
                                   {\normalfont\bf}}
\makeatother

\date{\today}

\newcommand{\iti}{\boldsymbol{1}}

\newcommand{\bbE}{\mathbb{E}}  
\newcommand{\bbP}{\mathbb{P}}  
\newcommand{\bR}{\mathbf{R}}

\newcommand{\dd}{\mathrm{d}}

\newcommand{\e}{\varepsilon}

\newtheoremstyle{new-thm}
 {3pt}
 {3pt}
 {\it}
 {0pt} 
 {\bf}
 {.}
 {.5em}
 {}
\newtheoremstyle{new-def}
 {3pt}
 {3pt}
 {\rm}
 {0pt} 
 {\bf}
 {.}
 {.5em}
 {}

\theoremstyle{new-thm}
  \newtheorem{thm}{Theorem}
  
  \newtheorem{lemma}[thm]{Lemma}

\theoremstyle{new-def}

  \newtheorem{rem}[thm]{Remark}

\begin{document}

\vspace*{2cm}
\begin{center} 
{\Large\bf Asymptotics of the probability distributions \\
of the first hitting times of Bessel processes}
\end{center} 

\bigskip

\begin{center} 
Yuji Hamana$^{\text{\rm a}}$ and 
Hiroyuki Matsumoto$^{\text{\rm b}}$ 
\footnote{Corresponding author.\\
{\it E-mail adress}: matsu(at)gem.aoyama.ac.jp(H.Matsumoto)} 
\end{center} 

\bigskip

\begin{center}
$^{\text{\rm a}}$Department of Mathematics, 
Kumamoto University, 
Kurokami 2-39-1, Kumamoto 860-8555, Japan \\
$^{\text{\rm b}}$Department of Physics and Mathematics, 
Aoyama Gakuin University, Fuchinobe 5-10-1, Sagamihara 252-5258, 
Japan 
\end{center}

\begin{quote} {\bf Abstract.} 
The asymptotic behavior of the tail probabilities 
for the first hitting times of the Bessel process 
with arbitrary index 
is shown without using the explicit expressions for the distribution 
function obtained in the authors' previous works.  

{\it keywords}: Bessel process, hitting time, tail probability
\end{quote}

\section{Introduction and main results}

Let $\bbP_a^{(\nu)}$ be the probability law 
on the path space $W=C([0,\infty);\mathbf{R})$ 
of a Bessel process with index $\nu\in\mathbf{R}$ 
or dimension $\delta=2(\nu+1)$ starting from $a>0$.   
For $w\in W$ we denote the first hitting time to $b>0$ 
by $\tau_b=\tau_b(w)$: 
\begin{equation*}
\tau_b=\inf\{t>0; w(t)=b\}.
\end{equation*}
\indent In recent works \cite{HM-I,HM-T} the authors have shown 
explicit forms of the distribution function and the density 
of the distribution of $\tau_b$ under $\bbP_a^{(\nu)}$ 
in the case where $0<b<a$.   
The other case, which is easier since we do not need to consider 
a natural boundary, has been known by Kent \cite{Kent}.   
See also Gettor-Sharpe \cite{GS}.  

When $\nu>0$ and $b<a$, it is shown in \cite{HM-T} that 
there exists a positive constant $C(\nu)$ such that 
\begin{equation*}
\bbP_a^{(\nu)}(\tau_b>t) = 1-\Bigl(\frac{b}{a}\Bigr)^{2\nu} + 
C(\nu) t^{-\nu} + o(t^{-\nu}).
\end{equation*}
The constant $C(\nu)$ may be expressed explicitly and 
we could treat all the cases.   
However, we need to consider separately the case 
where $\delta$ is an odd integer and 
the expression for $C(\nu)$ is different from the othere cases.  
This is because the expression for the distribution function 
itself is different.   

The aim of this note is to show that the constant $C(\nu)$ has 
the same simple expression also when $\delta$ is an odd integer 
by considering the asymptotics without using the explicit 
expressions for the distribution functions 
obtained in \cite{HM-T}.   

When $a<b$ and $\nu>0$, the explicit expressions for 
$\bbP_a^{(\nu)}(\tau_b>t)$ has been shown in \cite{Kent} 
and, from his result, it is easily shown that 
the tail probability decays exponentially.   
Hence we concentrate on the case of $b<a$.   

\begin{thm} \label{thm-1}
Let $\nu>0$ and $0<b<a.$   
Then{\rm,} as $t\to\infty,$ it holds that 
\begin{equation*}
\bbP_a^{(\nu)}(t<\tau_b<\infty)=
b^{2\nu}\Bigl\{1-\Bigl(\frac{b}{a}\Bigr)^{2\nu}\Bigr\}
\frac{1}{\Gamma(1+\nu)(2t)^{\nu}} 
+O(t^{-\nu-\e})
\end{equation*}
for any $\e\in(0,\frac{\nu}{1+\nu}),$ 
where $\Gamma$ denotes the usual Gamma function{\rm.}
\end{thm}

\begin{thm} \label{thm-2}
Let $\nu<0$ and $0<b<a.$   
Then{\rm,} as $t\to\infty,$ it holds that 
\begin{equation*}
\bbP_a^{(\nu)}(\tau_b>t)=
a^{2|\nu|}\Bigl\{1-\Bigl(\frac{b}{a}\Bigr)^{2|\nu|}\Bigr\}
\frac{1}{\Gamma(1+|\nu|)(2t)^{|\nu|}} 
+O(t^{-|\nu|-\e})
\end{equation*}
for any $\e\in(0,\frac{\nu}{1+\nu}).$
\end{thm}

When $\nu=0$, it is known that 
\begin{equation*}
\bbP_a^{(0)}(\tau_b>t)=\frac{2\log(a/b)}{\log t}+
o((\log t)^{-1}).
\end{equation*}
This identity has been discussed in \cite{HM-T} and 
we omit the details.   

\section{Proof of Theorem\ref{thm-1}}

We assume $\nu>0$ in this section.  
At first we give some lemmas.  
The first one is shown by Byczkowski and Ryznar \cite{BR}.

\begin{lemma} \label{l:rough-est}
There exists a constant $C$ such that 
\begin{equation*}
\bbP_a^{(\nu)}(t<\tau_b<\infty) \leqq Ct^{-\nu}.
\end{equation*}
\end{lemma}

\begin{lemma} \label{l:exp-tail}
If $0<b<a,$ one has 
\begin{equation} \label{e:exp-tail}
\bbP_a^{(\nu)}(\tau_b>t) = 1-\Bigl(\frac{b}{a}\Bigr)^{2\nu}+
\bbE_a^{(\nu)}\Bigl[ \Bigl(\frac{b}{R_t}\Bigr)^{2\nu}
\iti_{\{\inf_{0\leqq s \leqq t}R_s>b\}} \Bigr]
\end{equation}
for any $t>0,$ where $\bbE_a^{(\nu)}$ is the excpectation 
with respect to $\bbP_a^{(\nu)}$ and 
$\{R_s\}_{s\geqq0}$ denotes the coordinate process{\rm.}
\end{lemma}

\noindent{\bf Proof.}\ 
It is well known that 
\begin{equation*}
\bbP_a^{(\nu)}(\tau_b=\infty)=
\bbP_a^{(\nu)}(\inf_{s\geqq0}R_s>b)=
1-\Bigl(\frac{b}{a}\Bigr)^{2\nu}.
\end{equation*}
By the Markov property of Bessel processes, 
we have for $t>0$ 
\begin{align*}
\bbP_a^{(\nu)}(\tau_b=\infty) & = 
\bbP_a^{(\nu)}(\inf_{0\leqq s\leqq t}R_s>b \ \text{and}\ 
\inf_{s\geqq t}R_s>b) \\
 & = \bbE_a^{(\nu)}[ \bbP_{R_t}^{(\nu)}(\tau_b=\infty) 
\iti_{\{\inf_{0\leqq s \leqq t}R_s>b\}}] \\
 & = \bbE_a^{(\nu)}\Bigl[ 
\Bigl\{ 1-\Bigl(\frac{b}{R_t}\Bigr)^{2\nu} \Bigr\} 
\iti_{\{\inf_{0\leqq s \leqq t}R_s>b\}} \Bigr],
\end{align*}
which implies \eqref{e:exp-tail}. 

\begin{lemma} \label{l:asymp-mom}
For any $a>0$ and $p$ with $0<p<1+\nu,$ it holds that 
\begin{equation} \label{e:asymp-mom}
\frac{\Gamma(1+\nu-p)}{\Gamma(1+\nu)} \frac{1}{(2t)^{p}} 
e^{-\frac{a^2}{2t}} \leqq 
\bbE_a^{(\nu)}[(R_t)^{-2p}] \leqq 
\frac{\Gamma(1+\nu-p)}{\Gamma(1+\nu)} \frac{1}{(2t)^{p}} + 
Ct^{-1-p}
\end{equation}
for $t\geqq1,$ where $C$ is a positive constant 
independent of $t.$  
\end{lemma}

\noindent{\bf Proof.}\ 
By the explicit expression for the transition density 
of the Bessel process, we have 
\begin{equation*}
\bbE_a^{(\nu)}[(R_t)^{-2p}] = \int_0^\infty 
y^{-2p} \frac{1}{t} \Bigl(\frac{y}{a}\Bigr)^\nu y 
e^{-\frac{a^2+y^2}{2t}} I_\nu\Bigl(\frac{ay}{t}\Bigr) \dd y, 
\end{equation*}
where $I_\nu$ is the modified Bessel function of the first kind 
with index $\nu$ (cf \cite{W}) given by 
\begin{equation*}
I_\nu(z)=\Bigl(\frac{z}{2}\Bigr)^\nu \sum_{n=0}^\infty 
\frac{(z/2)^{2n}}{\Gamma(n+1) \Gamma(1+\nu+n)}
\quad (z\in\bR\setminus(-\infty,0)).
\end{equation*}
Hence, it is easy to get 
\begin{equation*}
\bbE_a^{(\nu)}[(R_t)^{-2p}] = \frac{1}{(2t)^{p}} e^{-a^2/2t} 
\sum_{n=0}^\infty \frac{a^{2n} \Gamma(n+\nu+1-p)}
{\Gamma(n+1)\Gamma(1+n+\nu)(2t)^n}
\end{equation*}
and the assertion of the lemma.

\begin{rem}
The moments of $R_t$ for fixed $t$ have 
explicit expressions by means of the Whittaker functions 
(cf. \cite{GR}, p.709), but it does not seem useful.
\end{rem}

We are now in a position to give a complete proof 
of Theorem \ref{thm-1}.   

\bigskip

\noindent{\bf Proof of Theorem \ref{thm-1}.}\ 
By Lemma \ref{l:exp-tail} we have 
\begin{equation*} 
\bbP_a^{(\nu)}(t<\tau_b<\infty)=
b^{2\nu}\bbE_a^{(\nu)}[(R_t)^{-2\nu}]-
b^{2\nu}\bbE_a^{(\nu)}[(R_t)^{-2\nu}\iti_{\{\tau_b\leqq t\}}].
\end{equation*}
For the first term we have shown in Lemma \ref{l:asymp-mom}
\begin{equation*} 
\bbE_a^{(\nu)}[(R_t)^{-2\nu}]=\frac{1}{\Gamma(1+\nu)(2t)^\nu}
(1+O(t^{-1})).
\end{equation*}
Hence, if we could show 
\begin{equation} \label{e:to-go}
\bbE_a^{(\nu)}[(R_t)^{-2\nu}\iti_{\{\tau_b\leqq t\}}]
=\frac{1}{\Gamma(\nu+1)(2t)^\nu} \Bigl(\frac{b}{a}\Bigr)^{2\nu} 
+O\Bigl(\frac{1}{t^{\nu+\e}}\Bigr) 
\end{equation}
for any $\e\in(0,\frac{\nu}{\nu+1})$, 
we obtain the assertion of the theorem.   

For this purpose, we let $\alpha\in(0,\frac{1}{\nu+1})$, 
choose $p$ satisfying 
\begin{equation*}
\frac{1}{1-\alpha}<p<\frac{1+\nu}{\nu}
\end{equation*}
and let $q$ be such that $p^{-1}+q^{-1}=1$.   
We devide the expectation on the right hand side of 
\eqref{e:to-go} into the sum of 
\begin{equation*}
I_1=\bbE_a^{(\nu)}[(R_t)^{-2\nu}
\iti_{\{\tau_b\leqq t^{\alpha q}\}}]
\quad \text{and} \quad
I_2=\bbE_a^{(\nu)}[(R_t)^{-2\nu}
\iti_{\{t^{\alpha q}<\tau_b\leqq t\}}]
\end{equation*}
\indent We simply apply the H\"older inequality to $I_2$.  
Then we get 
\begin{align*}
I_2 & \leqq \bbE_a^{(\nu)}[(R_t)^{-2\nu}
\iti_{\{t^{\alpha q}<\tau_b<\infty\}}] \\
 & \leqq \Bigl\{ \bbE_a^{(\nu)}[(R_t)^{-2\nu p}] \Bigr\}^{1/p}
\Bigl\{ \bbP_a^{(\nu)}(t^{\alpha q}<\tau_b<\infty) \Bigr\}^{1/q}.
\end{align*}
and, by Lemmas \ref{l:rough-est} and \ref{l:asymp-mom}, 
we see that there exists a constant $C_1$ such that 
\begin{equation*}
I_2 \leqq C_1 t^{-\nu-\alpha \nu}.
\end{equation*}
In the following we denote by $C_i$'s 
the constants independent of $t$.   

For $I_1$, the strong Markov property of Bessel processes implies
\begin{equation*}
I_1=\int_0^{t^{\alpha q}} \bbE_b^{(\nu)}[ (R_{t-s})^{-2\nu} ]
\bbP_a^{(\nu)}(\tau_b\in ds)=I_{11}+I_{12},
\end{equation*}
where 
\begin{equation*}
I_{11}=\int_0^{t^{\alpha q}} \frac{1}{2^\nu \Gamma(\nu+1)} 
\frac{1}{(t-s)^\nu} \bbP_a^{(\nu)}(\tau_b\in ds).
\end{equation*}
Since $\alpha q<1$, Lemma \ref{l:rough-est} imples 
\begin{equation*}
|I_{12}| = |I_1-I_{11}| \leqq \int_0^{t^{\alpha q}} 
\frac{C_2}{(t-s)^{\nu+1}} \bbP_a^{(\nu)}(\tau_b\in ds)
\leqq \frac{C_3}{t^{\nu+1}}.
\end{equation*}
\indent We devide $I_{11}$ into the sum of 
\begin{align*}
 & J_1=\int_{t^\alpha}^{t^{\alpha q}} \frac{1}{2^\nu\Gamma(\nu+1)} 
\frac{1}{(t-s)^\nu} \bbP_a^{(\nu)}(\tau_b\in ds) \\
\intertext{and} 
 & J_2=\int_{0}^{t^{\alpha}} \frac{1}{2^\nu\Gamma(\nu+1)} 
\frac{1}{(t-s)^\nu} \bbP_a^{(\nu)}(\tau_b\in ds). 
\end{align*}
For $J_1$ we have by Lemma \ref{l:rough-est} 
\begin{equation*}
0\leqq J_1 \leqq \frac{C_4}{(t-t^{\alpha q})^\nu} 
\bbP_a^{(\nu)}(t^\alpha<\tau_b<\infty) \leqq 
\frac{C_5}{t^{\nu+\alpha\nu}}.
\end{equation*}
\indent For $J_2$ we have 
\begin{align*}
J_2 & \leqq \frac{1}{2^\nu \Gamma(\nu+1) (t-t^\alpha)^\nu} 
\bbP_a^{(\nu)}(\tau_b\leqq t^\alpha) \\
 & \leqq \frac{1}{\Gamma(\nu+1)(2t)^\nu} 
\bbP_a^{(\nu)}(\tau_b<\infty)
\frac{1}{(1-t^{-(1-\alpha)})^\nu} \\
 & \leqq \frac{1}{\Gamma(\nu+1)(2t)^\nu}
\Bigl(\frac{b}{a}\Bigr)^{2\nu}
\Bigl( 1+\frac{C_6}{t^{1-\alpha}} \Bigr).
\end{align*}
On the other hand we have by Lemma \ref{l:rough-est}
\begin{align*}
J_2 & \geqq \frac{1}{\Gamma(\nu+1)(2t)^\nu} 
\bbP_a^{(\nu)}(\tau_b\leqq t^\alpha) \\
 & = \frac{1}{\Gamma(\nu+1)(2t)^\nu} \Bigl\{ 
\bbP_a^{(\nu)}(\tau_b<\infty) - 
\bbP_a^{(\nu)}(t^\alpha \leqq \tau_b<\infty) \Bigr\} \\
 & \geqq \frac{1}{\Gamma(\nu+1)(2t)^\nu} 
\Bigl(\frac{b}{a}\Bigr)^{2\nu} - \frac{C_7}{t^\nu t^{\alpha\nu}}.
\end{align*}
\indent Combining the above estimates, we obtain 
\begin{equation*}
\bbE_a^{(\nu)}[(R_t)^{-2\nu} \iti_{\{\tau_b\leqq t\}}] 
= \frac{1}{\Gamma(\nu+1)(2t)^\nu} \Bigl(\frac{b}{a}\Bigr)^{2\nu}
+ \frac{1}{t^{\nu}} O\Bigl(\frac{1}{t^{\alpha\nu}}+
\frac{1}{t}+\frac{1}{t^{1-\alpha}}\Bigr).
\end{equation*}
Since 
\begin{equation*}
0<\alpha\nu<\frac{\nu}{\nu+1}<1-\alpha<1
\end{equation*}
and we can choose arbitrary $\alpha$ satifying this condition, 
\begin{equation*}
\bbE_a^{(\nu)}[(R_t)^{-2\nu} \iti_{\{\tau_b\leqq t\}}] 
= \frac{1}{\Gamma(\nu+1)(2t)^\nu} \Bigl(\frac{b}{a}\Bigr)^{2\nu}
+ O\Bigl(\frac{1}{t^{\nu+\e}}\Bigr)
\end{equation*}
holds for any $\e\in(0,\frac{\nu}{\nu+1})$.  

Now we have shown \eqref{e:to-go} and 
the assertion of Theorem \ref{thm-1}.

\section{Proof of Theorem \ref{thm-2}}

Theorem \ref{thm-2} is easily obtained from Theorem \ref{thm-1}.   

We recall explcit expressions for 
the Laplace transforms of the distributions of $\tau_b$: 
for $\nu\in\mathbf{R}$, it is known (\cite{GS, HM-T}) that 
\begin{equation*}
\bbE_a^{(\nu)}[ e^{-\lambda \tau_b} ] = 
\Bigl(\frac{b}{a}\Bigr)^\nu 
\frac{K_\nu(a\sqrt{2\lambda})}{K_\nu(b\sqrt{2\lambda})}, 
\qquad \lambda>0, 
\end{equation*}
where $K_\nu$ is the modified Bessel function of the second kind. 
From this identity we easily obtain for $\nu>0$ 
\begin{equation*}
\bbP_a^{(-\nu)}(\tau_b\in dt) = \Bigl(\frac{a}{b}\Bigr)^{2\nu} 
\bbP_a^{(\nu)}(\tau_b \in dt).
\end{equation*}
Hence we get from Theorem \ref{thm-1} 
\begin{equation*} \begin{split}
\bbP_a^{(-\nu)}(\tau_b>t) & = \Bigl(\frac{a}{b}\Bigr)^{2\nu} 
\bbP_a^{(\nu)}(t<\tau_b<\infty) \\
 & = a^{2\nu} \Bigl\{ 1-\Bigl(\frac{b}{a}\Bigr)^{2\nu} 
\Bigl\} \frac{1}{\Gamma(1+\nu)(2t)^{\nu}} (1+o(1)).
\end{split} \end{equation*}

\noindent{\bf Acknowledgement.}\  
We thank Professor Yuu Hariya for valuable discussions.

\end{document}